\documentclass[10pt,conference]{IEEEtran}
\usepackage{amssymb}
\usepackage{amsmath}
\usepackage{amsfonts}
\usepackage{verbatim}
\usepackage{color}
\usepackage{multicol}
\usepackage{geometry}

\newtheorem{theorem}{Theorem}

\renewcommand{\S}{\mathbb{S}}

\begin{document}

\title{Deterministic guarantees for $L^{1}$-{reconstruction}: 
\\
A large sieve approach with geometric flexibility}
\author{\IEEEauthorblockN{Lu\'{\i}s Daniel Abreu} 
\IEEEauthorblockA{
Acoustics Research Institute,} 
\IEEEauthorblockA{
Wohllebengasse 12-14, Vienna A-1040, Austria. }%
\IEEEauthorblockA{
Email: labreu@kfs.oeaw.ac.at} \and \IEEEauthorblockN{Michael Speckbacher} %
\IEEEauthorblockA{Universit{\'e} de Bordeaux,} 
\IEEEauthorblockA{
Cours de la Lib{\'e}ration,
F-33405 Talence, France.}%
\IEEEauthorblockA{
Email: michael.speckbacher@u-bordeaux.fr} }
\maketitle

\vspace{-1cm} 
{\footnotesize \textbf{\textit{Abstract -- } We present estimates of the $p$%
-concentration ratio for various function spaces on different geometries
including the line, the sphere, the plane, and the hyperbolic disc, using
large sieve methods. Thereby, we focus on $L^{1}$-estimates which can be
used to guarantee the reconstruction from corrupted or partial information.}}

\section{Introduction}

Consider a measure space $(X,\mu )$, a measurable subset $\Omega \subset X$
and a subspace of $\mathcal{S\subset }L^{p}(X,\mu )$. A fundamental quantity
in mathematical signal analysis is the $p$-concentration ratio, defined as 
\begin{equation}
\lambda _{p}(\Omega ,\mathcal{S}):=\sup_{f\in \mathcal{S}}\frac{\Vert f\cdot
\chi _{\Omega }\Vert _{p}^{p}}{\Vert f\Vert _{p}^{p}}.
\label{general-conc-probl}
\end{equation}

For $p=2$, the study of this quantity was the cornerstone of the body of
work nowadays referred to as the `Bell Lab papers' of Landau, Slepian and
Pollak, culminating on Landau's necessary conditions on sampling and
interpolation \cite{Land}. This was later extended to a variety of contexts,
including spaces of analytic functions \cite{Seip0}, and wavelet and Gabor
spaces \cite{Daubechies,DeMarie}.

For $p=1$, $\lambda _{1}(\Omega ,\mathcal{S})<\frac{1}{2}$ implies that%
\begin{equation}
\Vert f\cdot \chi _{\Omega }\Vert _{1}<\frac{1}{2}\Vert f\Vert _{1}\text{%
,}  \label{main}
\end{equation}%
meaning that \emph{every signal }$f\in S$\emph{\ is sparse (poorly
concentrated) in }$\Omega $. Under such conditions, a remarkable phenomenon
discovered by Logan \cite{Logan, DS}, holds: if we sense a signal $f$
corrupted by \emph{unknown} noise $N$ supported in an \emph{unknown} $\Omega 
$, then $f$ can be perfectly reconstructed as the solution of the $l_{1}$
minimization problem%
\begin{equation*}
f=\arg \min_{s\in S}\Vert (f+N)-s\Vert _{1}.
\end{equation*}%
Another reconstruction scenario can be derived by generalizing an
observation of Donoho and Stark for bandlimited discrete functions \cite[%
Theorem 9]{DS}. In the absence of noise, if we only sense the projection of
a general $f\in S$ on $\Omega ^{c}$, then $f$ can be perfectly reconstructed
as the solution of the $L^{1}$-minimization problem%
\begin{equation*}
f=\arg \min_{s\in S}\Vert s\Vert _{1}\text{, \ \ \ \ \ subject to }P_{\Omega
}s=P_{\Omega }f\text{.}
\end{equation*}%
The above scenarios offer the possibility of reconstructing a function from
highly incomplete information, at the cost of obtaining a constant such that 
$\Vert f\cdot \chi _{\Omega }\Vert _{1}<C\Vert f\Vert _{1}$, with $C\leq 
\frac{1}{2}$. This is in general extremely difficult and the only sharp
result in this direction is provided by Tao's uncertainty principle for
signals of prime lenght \cite{Tao}. In Section 7.3 of \cite{DS}, it has been
suggested that considering a random $\Omega $, one could improve the
estimate $\lambda _{1}(\Omega ,\mathcal{S)}$. The full potential of this
suggestion has been later explored in the groundbreaking papers of Donoho 
\cite{Donoho} and Cand\'{e}s, Romberg and Tao \cite{CRT}, where it has been
shown that selecting both $\Omega $ and $S$ randomly, one could obtain $%
C\leq \frac{1}{2}$ with high probability, under reasonable assumptions on
the measure of $\Omega $. This lead to an intense research activity on the
topic nowadays known as Compressive Sensing \cite{FR}.

We will pursue a different strategy for the estimation of the quantity $%
\lambda _{1}(\Omega ,\mathcal{S})$, suggested by a ingenious application of
the large sieve principle to sparse recovery problems by Donoho and Logan 
\cite{DL}. This work was inspired by the following variation of the
classical large sieve inequality, due to Bombieri and quoted by Montgomery
in \cite[p. 562]{Montgomery}. Let $\mu $ be\ a periodic measure on the
circle and $f(\alpha )=\sum_{k=m+1}^{m+n}a_{k}e^{2\pi ik\alpha }$ a
trigonometric polynomial. Then:%
\begin{equation}
\left\Vert f\right\Vert _{L^{2}[(0,1),\mu ]}^{2}\leq (n+\frac{2}{\delta }%
)\left( \sup_{\alpha }\int_{\alpha }^{\alpha +\delta }d\mu \right)
\left\Vert f\right\Vert _{L^{2}(0,1)}^{2}\text{.}  \label{Bombieri}
\end{equation}%
In \cite{DL}, $L^{1}$-versions of (\ref{Bombieri}) are obtained, assuring
that (\ref{main}) holds if the measure $\mu $ is sparse (low concentration
on sets $\left[ \alpha ,\alpha +\delta \right] $). This lead us to pursue
the program of obtaining large sieve inequalities of the Donoho-Logan type
for $1\leq p<\infty $, in several signal analytic settings where the
original large sieve methods lack key ingredients (for instance, Beurling's
theory of extremal functions). In this note we will outline the first
results of this program. For convenience of presentation, consider the
general setting of the first paragraph. Denote by $\mathcal{S}_{K}$ the
reproducing kernel subspace of $L^{p}(X,\mu )$ consisting of all functions $%
f $ such that 
\begin{equation}
f(x)=\int_{X}K(x,y)f(y)d\mu (y),\ \ \forall x\in X,
\label{reprod-kernel-form}
\end{equation}%
for some Hermitian kernel $K$. We assume that $X$ is a metric space, and
that the kernel satisfies 
\begin{equation*}
\sup_{y\in X}\int_{X}|K(x,y)|d\mu (x)<\infty .
\end{equation*}%
Then, for $1\leq p<\infty $, every $f\in \mathcal{S}_{K}$ satisfies the
following inequality: 
\begin{equation}
\Vert f\cdot \chi _{\Omega }\Vert _{p}^{p}\leq \sup_{y\in X}\int_{\Omega
}|K(x,y)|d\mu (x)\cdot \Vert f\Vert _{p}^{p}.  \label{res-RKHS}
\end{equation}%
This is a simple consequence of the reproducing kernel equation (\ref%
{reprod-kernel-form}), since 
\begin{align*}
\int_{\Omega }|F(x)|d\mu (x)& \leq \int_{\Omega }\int_{X}|K(x,y)f(y)|d\mu
(y)d\mu (x) \\
& \leq \sup_{y\in X}\int_{\Omega }|K(x,y)|d\mu (x)\cdot \Vert f\Vert _{1}.
\end{align*}%
The statement for general $1\leq p<\infty $ then follows from the
Riesz-Thorin theorem using the trivial observation that $\Vert f\cdot \chi
_{\Omega }\Vert _{\infty }\leq \Vert f\Vert _{\infty }$.

From inequality (\ref{res-RKHS}) we conclude that (\ref{main}) holds as long
we can assure that 
\begin{equation*}
\sup_{y\in X}\int_{\Omega }|K(x,y)|d\mu (x)<C\leq \frac{1}{2}\text{.}
\end{equation*}%
This is, of course, completely impossible to do in the absence of more
information. However, in several situations where the kernel $K(x,y)$ is
explicitly known, one can obtain large sieve type inequalities which can be
used to obtain useful estimates from (\ref{res-RKHS}).

We will see some examples in the sections below, where the estimates are
given in terms of a measure of the sparsity of the set $\Omega $ involving
the quantity 
\begin{equation*}
\sup_{z\in X}\frac{|\Omega \cap B_{\varrho }(z,R)|}{|B_{\varrho }(z,R)|}%
\text{,}
\end{equation*}%
with $B_{\varrho }(z,R)$ the ball centered at $z$, measured in the metric $%
\varrho $. Here, $\varrho $ is fine tuned to the underlying geometry of the
space $X$, and $R$ is chosen according to the space $\mathcal{S}_{K}$. We
will see examples involving line, planar euclidean, planar hyperbolic and
spherical geometries. Due to the applications in $l_{1}$-minimization, we
will mostly focus on the estimation of $\lambda _{1}(\Omega ,\mathcal{S})$,
but the methods we use work for general $1\leq p<\infty $, with the
exception of the spherical case, where the problem for $p=1$ is still open.
Nevertheless, we will also present a $p=2$ large sieve\ bound for finite
spherical harmonic expansions, which may be useful in different
applications, taking into account the recent applications of the $p=2$ large
sieve bounds in superresolution on the so called well-separated case \cite%
{Moitra,AB}.

The results are presented in the following sequence. For reference we start
with one of the Donoho-Logan's large sieve inequalities for band-limited
functions and then present the main results on the finite spherical harmonic
setting. We then move to phase-space contexts and outline the results for
Gabor spaces from \cite{abspe17-sampta, AS} and a new result for Bergman
spaces which can be translated to the setting of Cauchy wavelets.

\section{Concrete Large Sieve Inequalities}

\subsection{Donoho-Logan's Large Sieve for the Paley-Wiener Space}

Consider the Paley-Wiener space of band-limited functions 
\begin{equation*}
PW_{W}^{p}:=\Big\{f\in L^{p}(\mathbb{R}):\ \mbox{supp}(\hat{f})\subseteq
\lbrack -\pi W,\pi W]\Big\},
\end{equation*}%
In \cite{DL}, Donoho and Logan introduced the following notion of maximum
Nyquist density: 
\begin{equation*}
\rho _{\mathbb{R}}(\Omega ,W)=|W|\cdot \sup_{t\in {\mathbb{R}}}\big|\Omega
\cap \lbrack t,t+1/W]\big|,
\end{equation*}%
and obtained the large sieve inequality%
\begin{equation*}
\Vert f\cdot \chi _{\Omega }\Vert _{1}\leq \frac{\pi }{2}\cdot \rho _{{%
\mathbb{R}}}(\Omega ,W)\cdot \Vert f\Vert _{1}\text{.}
\end{equation*}%
This shows that $\rho _{\mathbb{R}}(\Omega ,W)<\frac{1}{\pi }$ is enough to
assure perfect recovery in the context outlined in the introduction. The
results in \cite{DL} also cover discrete settings and applications of the $%
p=2$ inequality.

\subsection{Finite Spherical Harmonics Expansions}

Let $\S ^2$ be the unit sphere in $\mathbb{R}^3$ and $\mathcal{S}_L$ be the
space of finite spherical harmonics expansions of maximum degree $L$, i.e.
if $Y_l^m$ denotes the spherical harmonics, then 
\begin{equation*}
\mathcal{S}_L:=\Big\{f:\S ^2\rightarrow\mathbb{C}:\
f=\sum_{l=0}^L\sum_{m=-l}^l a_l^m Y_l^m,\ a_l^m\in\mathbb{C}\Big\}.
\end{equation*}
Estimates for the $p$-concentration problem are of particular interest for
example in geo-sciences where measurements like satelite images, or weather
data are not available on the whole sphere. The Bell-Lab approach to
concentration in $\mathcal{S}_L$ has numerically been applied in \cite{Sim}.

The maximum Nyquist density on $\S ^2$, tailored to $\mathcal{S}_L$ is
defined as 
\begin{equation}  \label{def-max-nyqu}
\rho_{\S ^2}(\Omega,L)=\sup_{y\in\S ^2}\frac{|\Omega\cap \mathcal{C}%
_{t_{L,L}}(y)|}{|\mathcal{C}_{t_{L,L}}(y)|},
\end{equation}
where the area $|\cdot|$ is measured w.r.t. the shift invariant surface
measure, $t_{L,L}$ denotes the largest zero of the Legendre polynomial $P_L$%
, and $\mathcal{C}_{t_{L,L}}(y)$ denotes the spherical cap with angle $%
\arccos t_{L,L}$ centered at $y$. Note that $t_{L,L}$ is an increasing
sequence converging to $1$. 


In \cite{hryspe19} estimates for the $p$-concentration problem are given for 
$p=1,2$. In particular, the result in the Hilbertian case reads: 
\begin{equation}  \label{eq:intro-thm}
\lambda_2( \Omega,\mathcal{S}_L ) \leq B_{L}\cdot\rho_{\S ^2}(\Omega,L),
\end{equation}
where 
\begin{align*}
B_{L}:=(1-t_{L,L})\left(\int_{t_{L,L}}^1 P_L^2(t)dt\right)^{-1}, \ \ L =
1,2, \ldots,
\end{align*}
is optimal within the chosen approach. The sequence $B_L$ is convergent with
limit 
\begin{align*}
\lim_{L\rightarrow\infty} B_{L} = J_1(j_{0,1})^{-2} \approx 3.71038068570948,
\end{align*}
where $j_{\alpha,m}$ denotes the m-th positive zero of the Bessel function
of the first kind $J_\alpha$. %
\newline
In the case $p=1$, no rigorous proof is given in \cite{hryspe19}. The
estimate with respect to $\rho_{\mathcal{S}^2}$ is nevertheless transferred
to an equivalent collection of inequalities which are shown to be true in
the limit $l\rightarrow \infty$. For small $l$ the conditions are
numerically verified, which therefore gives strong evidence that the
following estimate holds 
\begin{equation}  \label{l1-intro}
\lambda_1(\Omega,\mathcal{S}_L) \leq (2A-1)^{-1}\cdot\rho_{\mathbb{S}%
^2}(\Omega,L),
\end{equation}
where $A\approx 0.680460162465512$ is a unique positive solution of the
equation 
\begin{align*}
j_{0,1}A -2\sqrt{A}\ J_1\big( j_{0,1}\sqrt{A}\big) = \frac12
j_{0,1}-J_1(j_{0,1}).
\end{align*}

\subsection{Large Sieve Principles for Gabor Spaces}

Let $z=(x,\omega )\subset \mathbb{R}^{2}$, and define the time frequency
shift $\pi $ as $\pi (z)f(t)=e^{2\pi i\omega t}f(t-x)$. The short-time
Fourier transform (STFT) is defined as 
\begin{equation*}
V_{g}f(z)=\langle f,\pi (z)g\rangle =\int_{{\mathbb{R}}}f(t)e^{-2\pi i\omega
t}g(t-x)dt.
\end{equation*}%
In the following, we restrict the choice of windows $g$ to the class of
Hermite functions $h_{r}$. We define the planar maximum Nyquist density as 
\begin{equation}
\rho _{\mathbb{R}^{2}}(\Omega ,R):=\sup_{z\in \mathbb{R}^{2}}|\Omega \cap
(z+D_{R})|\text{.}  \label{planarNyquist}
\end{equation}

The main result of \cite{abspe17-sampta,AS} (where it is actually proved for 
$1\leq p<\infty $) is the following:

\begin{theorem}
Let $\Omega \subset \mathbb{R}^{2}$ and $V_{h_{r}}f\in L^{1}(\mathbb{R}^{2})$%
. For every $0<R<\infty $, and every $r$, 
\begin{equation}\label{thm1}
\Vert V_{h_{r}}f\cdot \chi _{\Omega }\Vert _{1}\leq \frac{\rho _{\mathbb{R}%
^{2}}(\Omega ,R)}{C_{r}(R)}\Vert V_{h_{r}}f\Vert _{1}\text{,}\ 
\end{equation}%
with $C_{r}(R)=1-P_{r}(R)e^{-\pi R^{2}}$, and $P_{r}$ is a polynomial of
degree $2r$ satisfying $P_{r}(0)=1$.
\end{theorem}

The result crucially relies on the following local reproducing formula 
\begin{equation}
V_{h_{r}}f(z)=\frac{1}{C_{r}(R)}\int_{z+D_{R}}\hspace{-0.3cm}%
V_{h_{r}}f(w)\langle \pi (w)h_{r},\pi (z)h_{r}\rangle dw\text{,}
\label{double-orthog}
\end{equation}%
which is shown in \cite{AS} via the correspondence between the STFT with
Hermite windows and polyanalytic Bargmann-Fock spaces \cite{AbrGr}.

As an illustrative application of the above theorem in the case $r=0$, i.e. the case of Gaussian
window $\varphi=h_0$,  suppose that one observes only
the time-frequency content of a STFT outside a region $\Omega $, $%
H:=P_{\Omega ^{c}}V_{\varphi}f\in L^{1}(\mathbb{R}^{2})$, and that $\Omega $
satisfies $\rho _{\mathbb{R}^{2}}(\Omega ,R)<(1-e^{-\pi R^{2}})/2$. Then: 
\begin{equation*}
V_{\varphi }f=\underset{V_{\varphi}g\in L^{1}(\mathbb{R}^{2})}{\text{argmin}}%
\big\Vert V_\varphi g\big\Vert_{1}\text{, \ subject to }P_{\Omega
^{c}}\left( V_{\varphi}g\right) =H\text{.}
\end{equation*}

\subsection{The Hyperbolic Case: Bergman Spaces and Analytic Wavelets}

Let $\varrho _{\mathbb{D}}$ denote the pseudohyperbolic metric in the disc $%
\mathbb{D}$ 
\begin{equation*}
\varrho _{\mathbb{D}}(z_{1},z_{2})=\left\vert \frac{z_{1}-z_{2}}{1-z_{1}%
\overline{z_{2}}}\right\vert, 
\end{equation*}%
and let $B_{\varrho _{\mathbb{D}}}(z,R)$\ be the pseudohyperbolic ball of center 
$z\in \mathbb{D}$ and radius $R<1$ defined as $B_{\varrho _{\mathbb{D}%
}}(z,R)=\left\{ w\in \mathbb{D}:\varrho _{\mathbb{D}}(w,z)<R\right\} $. Moreover, we define the hyperbolic measure of a set $\Omega$ as
$$
|\Omega|_\mathbb{D}:=\int_\Omega (1-|z|^2)^{-2}dz.
$$ The
Bergman space $A_{\alpha }^{p}(\mathbb{D})$ \cite{He} on the unit disc is defined as
the space of all analytic functions $F$ on $\mathbb{D}$ such that 
\begin{equation*}
\lVert F\rVert _{A_{\alpha }^{p}(\mathbb{D})}^{p}=\frac{1}{\pi }\int_{%
\mathbb{D}}\left\vert f(z)\right\vert ^{p}(1-\left\vert z\right\vert
^{2})^{\alpha -2}\,dz<\infty.
\end{equation*}%
The reproducing
kernel of $A_{\alpha }^{2}(\mathbb{D})$ is given by
$$
\mathcal{K}_{\mathbb{D}}^{\alpha }(z,w)=(1-z\overline{w})^{-\alpha }.
$$
One can define a hyperbolic maximum Nyquist density as%
\begin{equation*}
\rho _{\mathbb{D}}(\Omega ,R):=\sup_{z\in \mathbb{D}}|\Omega \cap B_{\varrho
_{\mathbb{D}}}(z,R)|_\mathbb{D}\text{.}
\end{equation*}%

It is then possible to obtain a hyperbolic analogue of (\ref{thm1}).

\begin{theorem}
Let $\Omega \subset \mathbb{D}$ and $F\in A_{\alpha }^{1}(\mathbb{D})$. For
every $R<1$, 
\begin{equation}\label{A1-estimate}
\Vert F\cdot \chi _{\Omega }\Vert _{A_{\alpha }^{1}(\mathbb{D})}\leq 
\frac{2\rho _{\mathbb{D}}(\Omega ,R)}{C^{\alpha }(R)}\cdot\Vert F\Vert _{A_{\alpha }^{1}(%
\mathbb{D})}\text{,}\ 
\end{equation}
where $C^\alpha(R)=\frac{1}{\alpha-1}\left(1-(1-R^2)^{\alpha-1}\right).$
\end{theorem}

The following local reproducing formula obtained by Seip \cite[Theorem 2.6]%
{Seip0} plays a key role in the proof:%
\begin{equation}\label{LRHyperbolic}
f(z)=\frac{1}{C^{\alpha }(R)}\int_{B_{\varrho _{\mathbb{D}}}(z,R)}f(z)%
\mathcal{K}_{\mathbb{D}}^{\alpha }(z,w)(1-\left\vert w\right\vert
^{2})^{\alpha -2}\,dw\text{.}  
\end{equation}%
For $p=2$, the Bergman space $A^2_\alpha(\mathbb{D})$ is conformally equivalent to the Bergman
space on the upper half-plane $A_{\alpha }^{2}(\mathbb{C}^{+})$. The spaces $A_{\alpha }^{p}(\mathbb{C}^{+})$ can
be understood, up to a weight, as the phase space of a continuous \emph{%
wavelet} transform with analyzing wavelets of the form 
\begin{equation*}
\widehat{g_{\alpha }}(\xi ):=\frac{2^{(\alpha -1)/2}}{\Gamma (\alpha
-1)^{1/2}}\,\xi ^{(\alpha -1)/2}e^{-\xi }\chi _{\lbrack 0,\infty )}(\xi )%
\text{.}
\end{equation*}%
In that case, using the conformal map between $\mathbb{D}$ and $\mathbb{C}^{+}$, the reproducing formula %
\eqref{LRHyperbolic} can be moved to $A_{\alpha }^{1}(\mathbb{C}^{+})$ and an equivalent estimate as  \eqref{A1-estimate} can be shown.
Details will be given elsewhere, together with the extension to the class of
wavelets which have phase space representations in polyanalytic Bergman
spaces \cite{AJFAA}. 

\section*{Acknowledgement}

L.D. Abreu was supported by the Austrian Science Fund (FWF) project
'Operators and Time-Frequency Analysis' (P 31225-N32), and M. Speckbacher
was supported by an Erwin-Schr{\"{o}}dinger Fellowship (J-4254) of the FWF.


\end{document}